\documentclass[12pt]{birkmult}
\usepackage{amsfonts}
\usepackage{epsfig}
\usepackage{graphics}
\usepackage{amsmath,amsthm}
\usepackage{amssymb}
\usepackage{inputenc}
\usepackage[T2A]{fontenc}

\date{\today}

\renewcommand{\theequation}{\thesection.\arabic{equation}}
\newtheorem{theorem}{Theorem}[section]

\newtheorem{definition}[theorem]{Definition}

\newtheorem{remark}[theorem]{Remark}

\newtheorem{proposition}[theorem]{Proposition}

\newtheorem{problem}[theorem]{Problem}

\begin{document}

\title[Abstract Interpolation Problem]
{Abstract Interpolation Problem and Some Applications. II:
Coefficient Matrices.}

\author[Kheifets]{A. Kheifets}
\address{Department of Mathematical Sciences \\
University of Massachusetts Lowell \\
One University Avenue \\
Lowell, MA 01854, USA}
\email{Alexander\_Kheifets@uml.edu}

\dedicatory{Dedicated to Joe Ball on occasion of his 70-th birthday}



\begin{abstract}

The main content of this paper is Lectures 5 and 6 that continue lecture notes \cite{KhBer}.
Content of Lectures 1-4 of \cite{KhBer} is reviewed for the reader's convenience in sections 1-4, respectively.
It is shown in Lecture 5 how residual parts of the minimal unitary extensions, that correspond to solutions of the problem, yield some boundary properties of the coefficient matrix-function. These results generalize the classical Nevanlinna - Adamjan - Arov~- Krein theorem. Lecture 6 discusses how further properties of the coefficient matrices follow from denseness of certain sets in the associated function model spaces. The structure of the dense set reflects the structure of the problem data.

\end{abstract}

\maketitle

\section{Abstract Interpolation Problem.}

\subsection{Data of the Problem.}

Data of the Abstract Interpolation Problem consists of the following components:
a vector space $X$ (without any topology), a positive-semidefinite sesquilinear form $D$ on $X$,
linear operators $T_1$ and $T_2$ on $X$, separable Hilbert spaces $E_1$ and $E_2$,
linear mappings $M_1$ and $M_2$ from $X$ to $E_1$ and $E_2$ respectively. The data pieces are connected
by the following identity
\begin{equation}\label{1.10}
D(T_2x,T_2y)-D(T_1x,T_1y)=\langle M_1x,M_1y\rangle_{E_1}-\langle M_2x,M_2y\rangle _{E_2}.
\end{equation}

\noindent
Let $\mathbb D$ be the open unit disc, $|\zeta|<1$, and let $\mathbb T$ be the unit circle, $|\zeta|=1$.
Let $w(\zeta): E_1\to E_2$ be a contraction for every $\zeta\in\mathbb D$ and assume that $w(\zeta)$ is analytic
in variable $\zeta$. Functions of this type are called the Schur class (operator-valued) functions.

\subsection{de Branges-Rovnyak Function Space.}\label{S1.2}

Let $L^2(E_2\oplus E_1)$ be the space of vector functions on the unit circle $\mathbb T$ that are square summable
against the Lebesgue measure. Space $L^w$ is defined as the range of $L^2(E_2\oplus E_1)$ under
\begin{equation}\label{1.20}
\begin{bmatrix} I_{E_2} & w(t) \\ \\ w(t)^* & I_{E_1} \end{bmatrix}^{\frac{1}{2}}
\end{equation}
endowed with the range norm. de Branges - Rovnyak space $H^w$ is defined as a subspace of $L^w$
that consists of functions
$$
f=\begin{bmatrix} f_2 \\ \\ f_1 \end{bmatrix} \in L^w,\quad f_2\in H^2
_+(E_2),\quad f_1\in H^2_-(E_1).
$$

\subsection{Setting of the Problem.}

A Schur class function $w:E_1\to E_2$ is said to be a solution of the
AIP with data \eqref{1.10}, if there exists a linear mapping $F:X
\to H^w$ such that for all $x\in X$
\begin{eqnarray}
 i) && \Vert Fx\Vert^2_{H^w}\leq D(x,x);\label{1.30}\\
ii) && tFT_2x-FT_1x=\begin{bmatrix} I_{E_2} & w(t) \\ \\ w(t)^* & I_{E_1} \end{bmatrix}
\begin{bmatrix}-M_2x\\ \\ \ \ M_1x\end{bmatrix},\ {\rm a.e.} \ \
t\in{\mathbb T}.\label{1.40}
\end{eqnarray}
One can write $Fx$ as a vector of two components
$$
Fx=\begin{bmatrix} F_+ x\\ \\ F_-x\end{bmatrix}
$$
which are $E_2$ and $E_1$ valued, respectively. Then conditions $Fx\in H^w$ and $($i$)$
read as
\begin{eqnarray*}
(a)&& F_+x\in H^2_+(E_2),\\
(b)&& F_-x\in H^2_-(E_1),\\
(c)&& \Vert Fx\Vert^2_{L^w}\leq D(x,x).
\end{eqnarray*}
Sometimes we will call the pair $(w, F)$ a solution of the Abstract Interpolation Problem.

\subsection{Special Case.}

The following additional assumption on operators $T_1$ and $T_2$ is met in many concrete problems: the operators
\begin{equation}\label{1.50}
(\zeta T_2-T_1)^{-1} {\rm\ and\ } (T_2-\overline\zeta T_1)^{-1}
\end{equation}
exist for all $\zeta\in\mathbb D$ except for a discrete set. In this case condition $($ii$)$ can be written
as explicit formulae for $F_+$ and $F_-$
\begin{eqnarray}
(F^w_+x)(\zeta)&=&(w(\zeta)M_1-M_2)(\zeta T_2-T_1)^{-1}x \ ,\label{1.60}\\
(F_-^w x)(\zeta)&=&\overline\zeta(M_1-w(\zeta)^*M_2)(T_2-\overline\zeta T_1)^{-1}x \ . \label{1.70}
\end{eqnarray}
From here one can see that under assumptions \eqref{1.50}, for every solution $w$ there exists only one $F$
that satisfies $($ii$)$.

\medskip\noindent
References to this section are \cite{KKY, 6, 9, KhY}.

\section{Examples.}
\setcounter{equation}{0}

\subsection{The Nevanlinna-Pick Interpolation Problem.}\label{Example 1}

In this section we recall several classical problems of analysis that can be included in the AIP scheme.

\begin{problem}\label{2.10}
Let $\zeta_1,\dots, \zeta_n ,\dots$ be a finite or infinite sequence of points in the unit disk
$\mathbb D$; let $w_1,\dots,w_n ,\dots$ be a sequence of complex numbers.
One is interested in describing all the Schur class functions $w$ such that
\begin{equation}\label{2.25}
w(\zeta_k)=w_k .
\end{equation}
\end{problem}

\medskip\noindent
The well-known solvability criterion for this problem is:
\begin{equation}\label{2.50}
\begin{bmatrix}\frac{1-\overline w_kw_j}{1-\overline\zeta_k \zeta_j}\end{bmatrix}^n_{k,j=1}\geq 0,
\ \  \text{for\ all}\ \ n\ge 1 .
\end{equation}
We specify data of the Abstract Interpolation Problem \eqref{1.10} as follows:
the space $X$ consists of all sequences
\begin{equation}\label{2.100}
 x=\begin{bmatrix} x_1 \\ \vdots \\ x_n \\ \vdots \end{bmatrix}
 \end{equation}
that have a finite number of nonzero components;
$$ D(x,y)=\sum_{k,j}\overline y_k\ \frac{1-\overline w_k w_j}{1-\overline\zeta_k \zeta_j} \
x_j\  , \quad x,y\in X ; $$
\begin{equation}\label{2.200}
T_1= \begin{bmatrix}\zeta_1 & & &  \\ & \ddots & &
\\ & & \zeta_n & \\  & && \ddots\end{bmatrix},\quad T_2=I_X;
\end{equation}
$$
E_1=E_2=\mathbb C^1;\quad
M_1x=\sum_j x_j,\quad M_2x=\sum_j w_j x_j.
$$
The operators $T_1$ and $T_2$ meet the special case assumption \eqref{1.50}. Therefore, for every solution $w$ there is only one corresponding mapping $F^w$ that
can be written in form \eqref{1.60}, \eqref{1.70}. It can be further explicitly computed as
$$ (F_+^wx)(\zeta)= \sum_j\frac{w(\zeta)-w_j}{\zeta-\zeta_j}x_j, $$
$$
(F^w_-x)(\zeta)=\overline\zeta\sum_j\frac{1-\overline{w(\zeta)}w_j}{1-\overline \zeta\zeta_j}x_j .
$$
Since $w$ has non-tangential boundary values, the function $F^wx$ extends to the boundary of the unit disk. Since for $|t|=1$ we have
$\overline t=1/t$, $(F^wx)(t)$ further simplifies as follows:
$$
 (F^wx)(t)=\begin{bmatrix} 1 & w(t)\\ \\ w(t)^* & 1 \end{bmatrix}
\begin{bmatrix} -\sum\limits_j\frac{w_jx_j}{t-\zeta_j} \\ \\
\sum\limits_j\frac{x_j}{t-\zeta_j}\end{bmatrix},\quad |t|=1 . $$
\begin{theorem}\label{2.300}
The solution set of the Nevanlinna-Pick problem coincides with the solution set of the
Abstract Interpolation Problem with data \eqref{2.100} - \eqref{2.200}. Moreover, for data of this
type, inequality \eqref{1.30} turns into equality
$$
 \Vert F^w x\Vert^2_{H^w} = D(x,x)
$$
for every solution $w$ and for every $x\in X$.
\end{theorem}

\medskip\noindent
This example can be viewed as a special case of the one in the next subsection.

\subsection{The Sarason Problem.}\label{Example 2}

\begin{problem}\label{2.350}
Let $H^2_+$ be the Hardy space of the unit disk.
Let $\theta$ be an inner function, $K_\theta=H^2_+\ominus
\theta H^2_+ \  , \ T^*_\theta x
=P_+\overline t x\ (x \in K_\theta),\ \  W^* $ be a contractive operator on $K_
\theta$ that commutes with $T^*_\theta:W^* T^*_\theta
=T^*_\theta W^* $. Find all the
Schur class functions $w$ such that
$$ W^*  x =P_+\overline w x  \ . $$
\end{problem}

\medskip\noindent
We specify here the data of the Abstract Interpolation problem as follows: $X=K_\theta$;
$$
D(x,x)=\Vert x\Vert^2_{_{K_\theta}} - \Vert W^*x\Vert^2_{_{K_\theta}}, x\in X ;
$$
\begin{equation}\label{2.400}
T_1=I_{K_\theta}, \quad T_2=T^*_\theta ;
\end{equation}
$$ E_1=E_2=\mathbb C^1,\quad
M_1 x=(W^*x)(0) ,\quad M_2x=x(0),$$
 where the latter notation stands for the value of an $H^2_+$ function at $0$.

 \medskip\noindent
 The operators $T_1$ and $T_2$ meet the special case assumption \eqref{1.50}. Therefore, for every solution $w$ there is only one corresponding mapping $F^w$ that
 can be explicitly computed as
 $$
 F^w x=\begin{bmatrix}1 & w\\ \overline w & 1 \end{bmatrix}
 \begin{bmatrix} x  \\  -W^*x\end{bmatrix},
 $$
 (when variable is on $\mathbb T$).
\begin{theorem}\label{2.500}
The solution set of the Sarason Problem \ref{2.350} coincides with the solution set of the
Abstract Interpolation Problem with data \eqref{2.400}. Moreover, for data of this
type, inequality \eqref{1.30} turns into equality
$$
 \Vert F^w x\Vert^2_{H^w} = D(x,x)
$$
for every solution $w$ and for every $x\in X$.
\end{theorem}
\begin{remark}\label{2.600}{\rm
Let $\zeta_1,\dots, \zeta_n ,\dots$ be a finite or infinite sequence of points in the unit disk
$\mathbb D$; let $w_1,\dots,w_n ,\dots$ be a sequence of complex numbers such that \eqref{2.50} holds.
Let $\theta$ be the Blaschke product with zeros $\zeta_k$, if the latter satisfy the Blaschke condition
and $\theta=0$ otherwise. Note that
$$
\frac{1}{1-t\overline\zeta_k}\in K_\theta.
$$
We define
$$
W^*\frac{1}{1-t\overline\zeta_k}=\frac{\overline w_k}{1-t\overline\zeta_k}.
$$
$W^*$ extends by linearity to a dense set in $K_\theta$ and further, due to \eqref{2.50}, to a contraction
on the whole $K_\theta$. The set of solutions $w$ of the Sarason Problem \ref{2.350} with this $\theta$ and this $W$
coincides with the set of solutions of the Nevanlinna-Pick Problem \ref{2.10}. However, the data of the Abstract Interpolation
Problem in \eqref{2.400} differ from the ones in \eqref{2.100}-\eqref{2.200}. Moreover,
the coefficient matrices $S$ in the description formula \eqref{4.140} for the solution sets
are different and the associated universal colligations \eqref{4.10}-\eqref{4.30} $\mathcal A_0$
are non-equivalent.}
\end{remark}

\medskip\noindent
References to this problem are \cite{6, KhY-0}.

\subsection{The Boundary Interpolation Problem.}\label{Example 3}

\begin{definition}\label{2.700}
A Schur class function $w$ defined on the unit disk $\mathbb D$
is said to have {\rm an angular derivative in the sense of Carath\'eodory}
at a point $t_0\in\mathbb T$ if there exists a nontangential unimodular
limit
$$
w_0=\underset{\zeta\to t_0}{\lim} w(\zeta),\quad |w_0|=1,
$$
and there exists a nontangential limit
$$
w'_0=\underset{\zeta\to t_0}{\lim} \frac{w(\zeta)-w_0}{\zeta-t_0}.
$$
\end{definition}
\begin{theorem}[Carath\'eodory-Julia]\label{2.800}
A Schur class function $w(\zeta)$ has an angular derivative at
$t_0\in{\mathbb T}$ if and only if
\begin{equation}\label{2.850}
D_{w,t_0}\overset{\rm def}{=}
\liminf_{\zeta\to t_0}\frac {1-|w(\zeta)|^2}{1-|\zeta|^2} <\infty
\end{equation}
$($here $|\zeta|< 1,\ \zeta\to t_0$ in an arbitrary way$)$. In this case
$$
w'_0=D_{w,t_0}\cdot\frac{w_0}{t_0}
$$
and
$$
\frac {1-|w(\zeta)|^2}{1-|\zeta|^2} \to\ D_{w,t_0}
$$
as $\zeta$ goes to $t_0$ nontangentially. Note that $D_{w,t_0}\ge 0$ and $D_{w,t_0}=0$ if and only if
$w(\zeta)$ is a constant of modulus 1.
\end{theorem}
\begin{theorem}\label{2.875}
A Schur class function $w$ has an angular derivative
in the sense of Carath\'eodory at a point $t_0\in{\mathbb T}$ if and only if there exists
a unimodular constant $w_0$ such that
\begin{equation}\label{2.900}
\left|\frac{w(t)-w_0}{t-t_0}\right|^2 + \frac{1-|w(t)|^2}{|t-t_0|^2}\in L^1
\end{equation}
against the Lebesgue measure $m(dt)$ on $\mathbb T$. In this case
$$
\int\limits_{\mathbb T}
\left(\left|\frac{w(t)-w_0}{t-t_0}\right|^2 + \frac{1-|w(t)|^2}{|t-t_0|^2}\right) m(dt)
=D_{w,t_0},
$$
where $D_{w, t_0}$ is the same as in \eqref{2.850}.
In particular, \eqref{2.900} implies that
$$
\frac{w-w_0}{t-t_0} \in H^2_+.
$$
\end{theorem}
\begin{problem}\label{2.1000}
Let $t_0$ be a point on the unit circle $\mathbb T$, let $w_0$ be a complex number, $|w_0|=1$,
and let $0\leq D<\infty$ be a given nonnegative number. One wants to
describe all the Schur class functions $w$ such that
$$
w(\zeta)\to w_0\quad {\rm as}\quad \zeta \to t_0
$$
nontangentially, and
$$D_{w,t_0}\leq D.$$
\end{problem}

\medskip\noindent
We specify here the data of the Abstract Interpolation Problem as follows:
$$
X  = \mathbb C^1,\quad D(x,x)=\overline xDx,
$$
\begin{equation}\label{2.1100}
T_1x  = t_0x,\quad  T_2x=x,
\end{equation}
$$
E_1 = E_2=\mathbb C^1,
$$
$$
M_1x=x,\quad M_2x=w_0x .
$$

\medskip\noindent
 The operators $T_1$ and $T_2$ meet the special case assumption \eqref{1.50}. Therefore, for every solution $w$ there is only one corresponding mapping $F^w$ that
 can be explicitly computed as
 $$
 (F^w x)(t)=\begin{bmatrix}1 & w(t)\\ \\ \overline {w(t)} & 1 \end{bmatrix}
 \begin{bmatrix} -\frac{w_0}{t-t_0}  \\ \\ \ \ \frac{1}{t-t_0}\end{bmatrix}x,\quad |t|=1.
 $$
 Direct computation shows that
 $$
 \Vert F^wx\Vert^2_{H^w}=\overline x D_{w, t_0}x.
 $$
 \begin{theorem}\label{2.1200}
The solution set of the Boundary Interpolation Problem coincides with the solution set of the
Abstract Interpolation Problem with data specified in \eqref{2.1100}.
\end{theorem}

\medskip\noindent
There is indeed inequality in \eqref{1.30} for
some solutions $w$ of the Abstract Interpolation Problem with data specified in \eqref{2.1100}. We will discuss this in more
detail in Section 5.

\medskip\noindent
References to this problem are \cite{KhHamb, Sarason}. For higher order analogue of Theorem \ref{2.800} and related boundary interpolation Problem \ref{2.1000} see
\cite{BoKh1, BoKh2}.

\section{Solutions of the Abstract Interpolation Problem.}
\setcounter{equation}{0}

\subsection{Isometry Defined by the Data.}

We say that two vectors $x_1$ and $x_2$ in $X$ are $D$ equivalent if
\begin{equation}\label{3.5}
D(x_1,y)=D(x_2,y),\quad\forall y\in X.
\end{equation}
We consider the vector space of equivalence classes $\{[x],\ x\in X\}$. We
define inner product between two equivalence classes as
\begin{equation}\label{3.7}
\langle [x],[y]\rangle {\buildrel {\rm def}\over =}\ D(x,y) \ .
\end{equation}
After completion we get a Hilbert space that will be denoted by $H_0$. We rewrite identity \eqref{1.10} as
\begin{equation}\nonumber
D(T_2x,T_2y)+\langle M_1x,M_1y\rangle _{E_1}=D(T_1x,T_1y)+\langle M_2x,M_2y\rangle _{E_2},
\end{equation}
or, using definition \eqref{3.7}, as
\begin{equation}\label{3.10}
\langle[T_2x],[T_2y]\rangle+\langle M_1x,M_1y\rangle _{E_1}=\langle [T_1x],[T_1y]\rangle+\langle M_2x,M_2y\rangle _{E_2}.
\end{equation}
We set
\begin{equation}\label{3.20}
d_V\ {\buildrel {\rm def}\over =}\ {\rm Clos}\left\{\begin{bmatrix}[T_1x]
\\ \\  M_1x\end{bmatrix},\ x\in X\right\}\subseteq H_0\oplus E_1 \
\end{equation}
and
\begin{equation}\label{3.30}
\Delta_V {\buildrel {\rm def}\over =}\ {\rm Clos}\left\{\begin{bmatrix}[T_2x]
\\ \\ M_2x\end{bmatrix},\ x\in X\right\}\subseteq H_0\oplus E_2 .
\end{equation}
We define a mapping $V:d_V\to\Delta_V$ by the formula
\begin{equation}\label{3.40}
V:\begin{bmatrix}[T_1x]\\ \\  M_1x\end{bmatrix}
\ {\buildrel {\rm def}\over\to}\
\begin{bmatrix}[T_2x]\\ \\ M_2x\end{bmatrix} \ .
\end{equation}
In view of \eqref{3.10}, $V$ is an isometry.
\begin{remark}
An arbitrary isometry $V$ from $H_0\oplus E_1$ to $H_0\oplus E_2$ may appear here under appropriate choice of the data in
\eqref{1.10}.
\if{
: $X=H_0\oplus E_1$,
$D(x,x)$ is the natural square of the norm,
$$
T_1(h_0\oplus e_1)= (P_{H_0}P_{d_V}( h_0\oplus e_1))\oplus 0,
$$
$$
M_1(h_0\oplus e_1)= P_{E_1}P_{d_V}( h_0\oplus e_1)),
$$
$$
T_2(h_0\oplus e_1)= (P_{H_0}VP_{d_V}( h_0\oplus e_1))\oplus 0,
$$
$$
M_2(h_0\oplus e_1)= P_{E_2}VP_{d_V}( h_0\oplus e_1)).
$$
}\fi
\end{remark}

\if{
Let $w$ be a solution of the AIP, that is there exists a mapping $F:X\to
H^w$ that possesses properties \eqref{1.30} and \eqref{1.40}. Combining \eqref{1.40} and \eqref{3.7},
we get
$$ \Vert Fx\Vert^2_{H^w}\leq D(x,x)=\Vert[x]\Vert^2_{H_0} \ . $$
Hence, $Fx$ depends only on the equivalence class $[x]$ and the following mapping is well defined
\begin{equation}\label{3.50}
G[x]\ {\buildrel {\rm def}\over =}\ Fx \
\end{equation}
and
$$ \Vert G[x]\Vert^2_{H^w}=\Vert Fx\Vert^2_{H^w}\leq\Vert[x]\Vert^2_{H_0} \ , $$
i.e. $G$ is a contraction. By continuity $G$ extends to a contraction from $H_0$ to $H^w$.
\fi

\subsection{Unitary Colligations, Characteristic Functions, Fourier Representations.}

Let $H,\ E_1,\ E_2$ be separable Hilbert spaces. A unitary mapping $\mathcal A$ of
$H\oplus E_1$ onto $H\oplus E_2$
\begin{equation}\label{3.55}
\mathcal A :H\oplus E_1\to H\oplus E_2
\end{equation}
is said to be a {\it unitary colligation}.
The space $H$ is called the {\it state space} of the colligation, $E_1$ is
called the {\it input space}, and $E_2$ is called the {\it output space}.
Both $E_1$ and $E_2$ are called {\it exterior spaces}. Sometimes it is convenient to write the colligation $\mathcal A$
as a block matrix:
\begin{equation}\label{3.60}
\mathcal A=\begin{bmatrix} A & B  \\ C & D\end{bmatrix}:
\begin{bmatrix} H \\  E_1 \end{bmatrix}\to \begin{bmatrix} H \\  E_2 \end{bmatrix}.
\end{equation}
The {\it characteristic function} of the unitary colligation is defined as
\begin{equation}\label{3.70}
w(\zeta)=D+\zeta C(I_H -\zeta A)^{-1}B.
\end{equation}
It is well defined  on $\mathbb D$ analytic contractive operator-valued function from $E_1$ to $E_2$.
The {\sl Fourier representation} of the space $H$ associated with the colligation $\mathcal A$ is defined as
\begin{equation}\label{3.80}
(\mathcal G h)(\zeta)=\begin{bmatrix} (\mathcal G_+ h)(\zeta)\\ \\ (\mathcal G_- h)(\zeta)\end{bmatrix}
= \begin{bmatrix} C(I_H -\zeta A)^{-1} h \\ \\ \overline\zeta B^*(I_H-\overline\zeta A^*)^{-1} h\end{bmatrix},\quad
h\in H,\quad \zeta\in\mathbb D .
\end{equation}
$\mathcal G$ maps space $H$ onto the de Barnges-Rovnyak space $H^w$ associated with the characteristic function $w$
(see Section \ref{S1.2}).
\begin{definition}\label{3.85}
We define the {\rm residual subspace} $H_{\rm res}\subseteq H$ of the colligation $\mathcal A$   as the maximal subspace of $H$ that
reduces $\mathcal A$ $($that is invariant for
$\mathcal A$ and $\mathcal A^*)$. Equivalently $H_{\rm res}$ can be defined as the maximal subspace of $H$ that is invariant for $A$ and $A^*$,
and $C|_{H_{\rm res}}=0$, $B^*|_{H_{\rm res}}=0$.  The {\rm simple part} of the space $H$ is defined as $H_{\rm simp}=H\ominus H_{\rm res}$.
A unitary colligation $\mathcal A$ is said to be {\rm simple} with respect to the exterior spaces $E_1$ and $E_2$ if $H_{\rm res}$ is trivial.
\end{definition}
\noindent The following fact will be of crucial importance to us in Lecture 5.
\begin{theorem}\label{3.90} Let $\mathcal G$ be the mapping defined in \eqref{3.80}. Then
$\mathcal G$ maps $H_{\rm simp}$ onto $H^w$ unitarily and $\mathcal G$ vanishes on
$H_{\rm res}$.
\end{theorem}

\subsection{Unitary Extensions of the Isometry $V$ and Solutions of the Problem.}

\begin{definition}\label{3.95}
We say that a unitary colligation $\mathcal A$ of form \eqref{3.60} is a unitary extension of the
isometry $V$, defined in \eqref{3.40}, if $H_0\subseteq H$ and
\begin{equation}\label{3.100}
\mathcal A|d_V =V.
\end{equation}
\end{definition}

\medskip\noindent
An extension $\mathcal A$ of $V$ is said to be {\sl minimal} if it does not have a nontrivial residual subspace in $H\ominus H_0$.
Note that minimal extension $\mathcal A$ may have a residual subspace in $H$, though.
\begin{theorem}\label{3.110}
A Schur class function $w:E_1\to E_2$ and a mapping $F: X\to H^w$ solve Abstract Interpolation Problem \eqref{1.30} - \eqref{1.40} if and only if
there exists a unitary colligation $\mathcal A$ of form \eqref{3.60} that minimally extends the isometry $V$ $($defined in \eqref{3.40} from the problem data$)$ such that
$w$ is the characteristic function of $\mathcal A$
\begin{equation}\label{3.120}
w(\zeta)=D+\zeta C(I_H -\zeta A)^{-1}B
\end{equation}
and
\begin{equation}\label{3.130}
Fx=\mathcal G [x],\quad  \forall x\in X,
\end{equation}
where $[x]$ is the $D$-equivalence class of $x$ $($defined in \eqref{3.5}$)$ and $\mathcal G$ is the Fourier representation of the colligation $\mathcal A$ $($defined in
\eqref{3.80}$)$.
\end{theorem}

\medskip\noindent
References to this section are \cite{KKY, 6, 9, KhY}.

\section{Parametric Description of the Solutions of the Abstract Interpolation Problem.}
\setcounter{equation}{0}

\subsection{Structure of Minimal Unitary Extensions of $V$.}

Let $V$ be an isometric colligation:
$$ V:d_V\to \Delta_V \ , \ \ d_V\subseteq H_0\oplus E_1\ , \ \ \Delta_V\subseteq H_0
\oplus E_2 \ . $$
Let $\mathcal A$ be a minimal unitary extension of $V$:
$$\mathcal A:H\oplus E_1\to H\oplus E_2\ ,\quad H\supseteq H_0\ , \quad\mathcal A| {d_V}=V
\ . $$
Let $d^\perp_V$ and $\Delta^\perp_V$ be the orthogonal complements of $d_V$ in
$H_0\oplus E_1$ and $\Delta_V$ in $H_0\oplus E_2$, respectively. Let
\begin{equation}\label{4.2}
H_1=H\ominus H_0.
\end{equation}
Then the orthogonal complement of $d_V$ in $H\oplus E_1$ is $H_1\oplus
d_V^\perp$ and the orthogonal complement of  $\Delta_V$ in $H\oplus E_2$ is
$H_1\oplus \Delta_V^\perp$. Since $\mathcal A$ is a unitary operator and it maps $d_V$
onto $\Delta_V$ ($\mathcal A| {d_V}=V$), $\mathcal A$ has to map the orthogonal
complement onto the orthogonal complement, i.e. $H_1\oplus d_V^\perp$ onto
$H_1\oplus\Delta_V^\perp$. Denote the restriction of $\mathcal A$ onto $H_1\oplus d_V
^\perp$ by $\mathcal A_1$. Thus, $\mathcal A_1$ is a unitary colligation,
$$
\mathcal A_1:H_1\oplus d_V^\perp\to H_1\oplus\Delta_V^\perp.
$$
Since $\mathcal A$ is a minimal extension of $V$, $\mathcal A_1$
is a simple colligation with respect to $d^\perp_V$ and $\Delta^\perp_V$.
Thus, the parameter of a minimal unitary extension $\mathcal A$ is an arbitrary simple unitary colligation $\mathcal A_1$
with input space $d^\perp_V$ and output space $\Delta^\perp_V$.

\medskip\noindent
Let $N_1$ be an auxiliary copy of $d_V^\perp$, that is we assume that there exists a unitary mapping $u_1$ from $d_V^\perp$ onto $N_1$
\begin{equation}\label{4.4}
u_1: d_V^\perp\to N_1.
\end{equation}
Also let
$N_2$ be a copy of $\Delta_V^\perp$, that is there exists a unitary mapping $u_2$ from $\Delta_V^\perp$ onto $N_2$
\begin{equation}\label{4.6}
u_2: \Delta_V^\perp\to N_2.
\end{equation}
In what follows it will be convenient to consider simple unitary colligations $\mathcal A_1$ with input space $N_1$ and output space $N_2$
(instead of $d_V^\perp$ and $\Delta_V^\perp$)
\begin{equation}\label{4.8}
\mathcal A_1:H_1\oplus N_1\to H_1\oplus N_2.
\end{equation}

\subsection{Universal Extension of $V$.}

Here we define a unitary colligation $\mathcal A_0$ that extends $V$ in a different way:
\begin{equation}\label{4.10}
\mathcal A_0:H_0\oplus E_1\oplus N_2\to H_0\oplus E_2\oplus N_1
\end{equation}
with
\begin{equation}\label{4.20}
\mathcal A_0| {d_V}=V, \quad (d_V\subseteq H_0\oplus E_1),
\end{equation}

\begin{equation}\label{4.30}
\mathcal A_0| {d_V^\perp}= u_1, \quad \mathcal A_0| {N_2}=u^*_2 ,
\end{equation}
where $u_1$ maps unitarily $d_V^\perp$ onto $N_1$ and $u_2$ maps unitarily $\Delta_V^\perp$ onto $N_2$.
$\mathcal A_0$ is uniquely defined by $V$ and the identification maps $u_1$ and $u_2$. Note that
mappings $u_1$ and $u_2$ can be chosen arbitrarily.  We will call this choice a {\sl normalization} of the
universal colligation $\mathcal A_0$.

\medskip\noindent
Similarly to \eqref{3.60}, we write colligation $\mathcal A_0$
as a block matrix:
\begin{equation}\label{4.60}
\mathcal A_0=\begin{bmatrix} A_0 & B_0  \\ C_0 & D_0\end{bmatrix}:
\begin{bmatrix} H_0 \\  E_1\oplus N_2 \end{bmatrix}\to \begin{bmatrix} H_0 \\  E_2\oplus N_1 \end{bmatrix}.
\end{equation}
We also introduce the characteristic function of $\mathcal A_0$
\begin{equation}\label{4.70}
S(\zeta)=D_0+\zeta C_0(I_{H_0} -\zeta A_0)^{-1}B_0.
\end{equation}
It is an analytic on $\mathbb D$ contractive operator-valued function from $E_1\oplus N_2$ to $E_2\oplus N_1$.
Note that $S$ depends on the data of the problem and on normalization \eqref{4.30} of $\mathcal A_0$.
According to the structure of the input an the output spaces, we further break $S$ into blocks
\begin{equation}\label{4.75}
S(\zeta)=\begin{bmatrix} s_0(\zeta) & s_2(\zeta)  \\ \\ s_1(\zeta) & s(\zeta)\end{bmatrix}:
\begin{bmatrix} E_1\\ \\ N_2 \end{bmatrix}\to \begin{bmatrix} E_2\\ \\ N_1 \end{bmatrix}.
\end{equation}
The special structure \eqref{4.20}, \eqref{4.30} of $\mathcal A_0$ forces
\begin{equation}\label{4.77}
s(0)=0.
\end{equation}
We consider the {\sl Fourier representation} of the space $H_0$ associated with the colligation $\mathcal A_0$
\begin{equation}\label{4.80}
(\mathcal G_0 h_0)(\zeta)=\begin{bmatrix} (\mathcal G_{0+} h_0)(\zeta)\\ \\ (\mathcal G_{0-} h_0)(\zeta)\end{bmatrix}
= \begin{bmatrix} C_0(I_{H_0} -\zeta A_0)^{-1} h_0 \\ \\ \overline\zeta B_0^*(I_{H_0}-\overline\zeta A_0^*)^{-1} h_0\end{bmatrix},
\end{equation}
$h_0\in H_0, \zeta\in\mathbb D.$
$\mathcal G_0$ maps space $H_0$ onto the de Barnges-Rovnyak space $H^S$ associated with the characteristic function $S$.

\subsection{Description of Solutions.}

Given unitary colligation $\mathcal A_0$ of form \eqref{4.10} and unitary colligation $\mathcal A_1$ of form \eqref{4.8}, there is a procedure (called
{\sl feedback coupling}) that produces a colligation $\mathcal A$ of the form \eqref{3.55} with
$$
H=H_0\oplus H_1.
$$
We do not discuss here the feedback coupling. A detailed explanation of the procedure is given, for instance, in \cite{KhBer}, Lecture 4, Section 3, page 373.
Here we just state some consequences of this procedure.
\begin{theorem}\label{4.85}
If $\mathcal A_0$ is defined in terms of $V$ as in \eqref{4.20}-\eqref{4.30} and $\mathcal A_1$ is an
arbitrary simple unitary colligation of form \eqref{4.8},
then the feedback coupling $\mathcal A$ is a
minimal unitary extension of $V$ in the sense of Definition \ref{3.95}. Moreover, all  minimal unitary extensions of $V$ arise in this way.
\end{theorem}
\begin{theorem}\label{4.90}
Given unitary colligation $\mathcal A_0$ of the form \eqref{4.10} and unitary colligation $\mathcal A_1$ of the form \eqref{4.8}.  Let unitary colligation $\mathcal A$
 of form \eqref{3.55} be the feedback coupling of $\mathcal A_0$ and $\mathcal A_1$. Then
\begin{equation}\label{4.100}
w=s_0+s_2\omega(I_{N_1}-s\omega)^{-1}s_1,
\end{equation}
where $w$ is the characteristic function of $\mathcal A$, $\omega$ is the characteristic function of $\mathcal A_1$, and $S$  is
the characteristic function of $\mathcal A_0$ $($see \eqref{4.75}$)$;
\begin{equation}\label{4.110}
\mathcal G\begin{bmatrix} h_1 \\ h_0\end{bmatrix}=\begin{bmatrix}\psi\omega &
I_{E_2} & 0 & 0 \\ 0 & 0 & \varphi^*\omega^* & I_{E_1}\end{bmatrix}
\mathcal G_0 h_0+\begin{bmatrix}\psi & 0 \\ 0 & \varphi^* \end{bmatrix}
\mathcal G_1 h_1 \ ,
\end{equation}
$ h_0\in H_0, h_1\in H_1$, where $\mathcal G$, $\mathcal G_1$ and $\mathcal G_0$ are Fourier representations of
$\mathcal A$, $\mathcal A_1$ and $\mathcal A_0$, respectively, as defined in \eqref{3.80}, \eqref{4.80},
\begin{equation}\label{4.120}
\varphi=(I_{N_1}-s\omega)^{-1}s_1,\quad
\psi=s_2(I_{N_2}-\omega s)^{-1}.
\end{equation}
\end{theorem}
\noindent
Combining Theorems \ref{3.110} and  \ref{4.90}, we get a description of all solutions
$(w,F)$ of the Abstract Interpolation Problem \eqref{1.30}-\eqref{1.40}.
\begin{theorem}\label{4.130}
Let $V$ be the isometry defined by the data of the problem as in \eqref{3.20} - \eqref{3.40}. Let $N_1$ and $N_2$ be the spaces \eqref{4.4}, \eqref{4.6} and
let $\mathcal A_0$ be the unitary colligation defined in \eqref{4.10} - \eqref{4.30}. Let $S$ be the characteristic function \eqref{4.75} of $\mathcal A_0$  and
$\mathcal G_0$ be the Fourier representation \eqref{4.80} of $H_0$. Then the
solution set $(w,F)$ of the Abstract Interpolation Problem \eqref{1.30}-\eqref{1.40} is described as follows
\begin{equation}\label{4.140}
w=s_0+s_2\omega(I_{N_1}-s\omega)^{-1}s_1,
\end{equation}
where $\omega$ is an arbitrary Schur class function from $N_1$ to $N_2$, $S$ is the characteristic function \eqref{4.75} of $\mathcal A_0$;
\begin{equation}\label{4.150}
Fx=\begin{bmatrix}\psi\omega &
I_{E_2} & 0 & 0 \\ 0 & 0 & \varphi^*\omega^* & I_{E_1}\end{bmatrix}
\mathcal G_0 [x],\quad x\in X,
\end{equation}
where $[x]$ is the $D$-equivalence class of $x$ defined in \eqref{3.5},
\begin{equation}\label{4.160}
\varphi=(I_{N_1}-s\omega)^{-1}s_1,\quad
\psi=s_2(I_{N_2}-\omega s)^{-1},
\end{equation}
$\omega$, $s$, $s_1$, $s_2$ are the same as above. All functions in formula \eqref{4.150} are considered on
$\mathbb T$.
\end{theorem}

\medskip\noindent
References to this section are \cite{ArGr, KKY, 6, 9, KhY}.

\section{Lecture 5: Inequality $\Vert Fx\Vert^2\leq D(x,x)$, Residual Parts of Minimal Unitary Extensions and the
Nevanlinna - Adamjan - Arov - Krein Type Theorems.}
\setcounter{equation}{0}

This lecture is focused on the inequality $\Vert Fx\Vert^2\leq
D(x,x)$ in the seting of the Abstract Interpolation  Problem (AIP)
\eqref{1.30}, \eqref{1.40}. The main goals are

\begin{itemize}
\item to explain the inequality in terms of the corresponding minimal unitary
extension \eqref{3.100} $\mathcal A$ of the isometry $V$ \eqref{3.20}-\eqref{3.40}. For more details, see \cite{6,9};

\medskip
\item to give a formula for the quantitative characteristic of how
the inequality is far from the equality.
\end{itemize}

\medskip\noindent
After that we apply the latter formula to the case when the equality $\Vert Fx
\Vert^2=D(x,x)$ is known a priori (like in Problem \ref{2.10} and Problem \ref{2.350}).
This in turn yields certain boundary properties
of the coefficient matrix $S(\zeta)$, defined in \eqref{4.75}.
In particular, this leads to generalizations of a classical Nevanlinna~- Adamjan~- Arov~- Krein theorem \cite{23, 24}: for general semi-determinate Nehari problem  \cite{9, 10, KhNeh} and for general Commutant Lifting problem  \cite{BallK1, BallK2}.

\subsection{The Equality $\Vert Fx\Vert^2=D(x,x)$ and Simplicity of the Corresponding Minimal Unitary Extension.}

\medskip
Let a Schur class operator function $w:E_1\to E_2$ and a mapping $F:X\to H^w$ be a solution of the AIP \eqref{1.30}, \eqref{1.40}.
\if{
with data \eqref{1.10}: linear space $X$, quadratic
form $D$, separable Hilbert spaces $E_1$ and $E_2$, linear operators
$T_1:X\to X, \ T_2:X\to X, \ M_1:X\to E_1,\ M_2:X\to E_2$. All these objects are
bound to one another by the Fundamental Identity.
}\fi
Let $V : d_V\to\Delta_V,\ d_V\subseteq H_0\oplus E_1,\ \Delta_V\subseteq H_0\oplus E_2$ be
the isometry \eqref{3.20}-\eqref{3.40} associated to AIP data \eqref{1.10}. Then, by Theorem \ref{3.100},
$w$ is the characteristic function of a unitary colligation $\mathcal A$ of the form \eqref{3.60} that minimally
extends the isometry $V$ and
\begin{equation}\label{5.100}
Fx=\mathcal G [x],\quad  \forall x\in X,
\end{equation}
where $[x]$ is the $D$-equivalence class of $x$ $($defined in \eqref{3.5}$)$ and $\mathcal G$ is the Fourier representation of the colligation $\mathcal A$ $($defined in
\eqref{3.80}$)$. By Theorem \ref{3.90},
$\mathcal G $ maps $H_{\rm simp}\subseteq H$ onto
$H^w$ unitarily and $\mathcal G $ vanishes on $H_{\rm res}\subseteq H$, see also Definition \ref{3.85}.

\medskip\noindent
From here one can see what the equality
\begin{equation}\label{5.200}
\Vert Fx\Vert^2_{H^w}=D(x,x), \ \forall x\in X
\end{equation}
means: in view of \eqref{5.100} and definition \eqref{3.5}-\eqref{3.7} of $H_0$, equality \eqref{5.200} is the same as
$$ \Vert\mathcal G [x]\Vert^2_{H^w}=\Vert[x]\Vert^2_{H_0}\ . $$
Since the lineal $\{[x],\ x\in X\}$ is dense in $H_0$, this means
that $\mathcal G $ is isometric on $H_0$, i.e., $H_0\subseteq H_{\rm simp}$. The latter
is equivalent to the inclusion $H_{\rm res}\subseteq H\ominus H_0$. Since
extension $\mathcal A$ is minimal, this is possible if and only if
$H_{\rm res}=\{0\}$. Thus, we arrive at the following
\begin{proposition}\label{5.300}
$\Vert Fx\Vert^2_{H^w}=D(x,x), \ \forall x\in X$ if and only if the
corresponding minimal unitary extension $\mathcal A$ of the isometry $V$ is simple.
\end{proposition}

\noindent
Hence, a strict inequality in \eqref{5.200} may occur for some $x\in X$ if and only if the corresponding
minimal extension $\mathcal A$ is non-simple, i.e., $H\supseteq H_{\rm res}\neq\{0\}$. In
this case $\mathcal A|H_{\rm res}$ is a unitary operator on $H_{\rm res}$.

\subsection{Residual Part of a Minimal Unitary Extension and its
Spectral Function.}

\medskip
The theme of this section is the unitary operator $\mathcal A|_{H_{\rm res}}$. By Theorem \ref{4.85},
every minimal unitary extension $\mathcal A:H\oplus
E_1\to H\oplus E_2$ of the isometry $V$ is the feedback
coupling of the universal unitary colligation $\mathcal A_0$, defined in \eqref{4.20}-\eqref{4.30}, and a simple unitary colligation $\mathcal A_1$
of form \eqref{4.8}:
$$ \mathcal A_0:H_0\oplus E_1\oplus N_2\to H_0\oplus E_2\oplus N_1, $$
$$\mathcal  A_1:H_1\oplus N_1\to H_1\oplus N_2, $$
where $H=H_0\oplus H_1$.
It follows from the feedback coupling procedure that if colligations $\mathcal A_0$ and $\mathcal A_1$ are not simple,
then their residual parts are contained in the residual part of their coupling $\mathcal A$. However, the residual part
of the latter colligation may be properly larger.

\medskip\noindent
Here we are interested in the piece of the residual part of $\mathcal A$ that
results from the feedback coupling procedure, but not from non-simplicity
of the coupled colligations. Let $\mathcal A_0:H_0\oplus E_1\oplus N_2\to H_0\oplus E_2\oplus N_1$
and $\mathcal A_1:H_1\oplus N_1\to H_1\oplus N_2$ be simple unitary colligations.
Let $\mathcal A:H\oplus E_1\to H\oplus E_2$ be their feedback coupling, where
$H=H_0\oplus H_1$. Let $H=H_{\rm simp}\oplus H_{\rm res}$ be the simple and
residual parts of the colligation $\mathcal A$, respectively. Let $U=\mathcal A|H_{\rm res}$.

\begin{definition}\label{5.400}
Let $U$ be a unitary operator on a separable
Hilbert space $\mathcal H$. Let $N$ be an auxiliary Hilbert space and let
$\Gamma :N\to \mathcal H $ be a linear operator from $N$ into $\mathcal H $ such that $\Gamma
(N)$ is a {\it cyclic subspace} for the operator $U$, i.e., the closed
linear span of $\{U^k\Gamma (N)\}_{k\in\mathbb Z}$ coincides with $\mathcal H $.
The operator function $a(\zeta):N\to N$, 
$$ a(\zeta)\ {\buildrel def\over = }\ \frac{1}{2} \Gamma ^*\ {{\bf 1}_{\mathcal H }+\zeta U
\over {\bf 1}_{\mathcal H }-\zeta U}\ \Gamma  $$
is called a {\it spectral function} of the operator $U$. Clearly, $a$ is analytic in $\mathbb D$.
\end{definition}

\medskip\noindent
Since $U$ is a unitary operator, $a(\zeta)+a(\zeta)^*\geq 0$. Therefore, $a(\zeta)$
admits a Riesz-Herglotz representation
$$ a(\zeta)=\frac{1}{2}\int\limits_{\mathbb T }{t+\zeta\over t-\zeta}\ \sigma (dt), $$
where $\sigma (dt)$ is an operator valued measure $(\sigma (dt):N\to N$) on the
unit circle $\mathbb T $. The measure $\sigma (dt)$ is called a {\it spectral
measure} of the operator $U$.

\begin{remark}\label{5.500}
Different choices of the auxiliary space $N$ and the operator $\Gamma $ lead to different
spectral functions. However, any of them defines the unitary operator $U$
uniquely up to a unitary equivalence.
\end{remark}

\medskip\noindent
The next theorem gives a formula for the spectral function of the unitary
operator $U=\mathcal A| H_{\rm res}$ that is the residual part of the feedback
coupling considered above in this section. Results of this type in the context of the
cascade coupling  go back to Yu. L. \v Smulian \cite{22}
and were inspired by M. Liv\v sits and M. Brodskii \cite{BrLiv}. It was realized,
in particular, as an obstacle for solving the Hilbert space
invariant subspace problem by means of factorization of the characteristic 
function.  In the same context it was investigated  by B. Sz.-Nagy and
C. Foia\c s (\cite{NF}, irreducible factorizations), and by L. de Branges (\cite{deBr},
overlapping subspaces), see also \cite{Schv, Ball1978}.

\begin{theorem}[ \cite{6, 9} ]\label{5.600}
Let unitary colligation $\mathcal A:H\oplus E_1
\to H\oplus E_2$ be the feedback coupling of simple unitary
colligations $\mathcal A_0:H_0\oplus E_1\oplus N_2\to H_0\oplus E_2\oplus N_1$
and $\mathcal A_1:H_1\oplus N_1\to H_1\oplus N_2$, where $H=H_0\oplus H_1$. Let
$$ 
S(\zeta)=\begin{bmatrix} s_0(\zeta) &  s_2(\zeta) \\ s_1(\zeta) &  s(\zeta)
\end{bmatrix}:\
\begin{bmatrix}E_1\cr N_2 \cr\end{bmatrix}\to\begin{bmatrix}
E_2 \cr N_1\cr\end{bmatrix} 
$$
be the characteristic function of the colligation $\mathcal A_0$ and $\omega (\zeta)$ be the
characteristic function of the colligation $\mathcal A_1$.
Let $H_{\rm res}\subseteq H$ be the
residual subspace of $\mathcal A$. Let 
$$
U=\mathcal A|_{H_{\rm res}}
$$
and let $N\ {\buildrel def\over
 = }\ N_2\oplus N_1$. We define $\Gamma :N\to H_{\rm res}$ as
$$ \Gamma \left(\begin{bmatrix}n_2 \\ n_1\cr\end{bmatrix}\right) \ {\buildrel
def\over = }\ P_{H_{\rm res}}(P_{H_0}\mathcal A_0^*(0_{H_0}\oplus 0_{E_2}\oplus(
\omega (0)^*n_2+n_1))\oplus P_{H_1}\mathcal A^*_1(0_{H_1}\oplus n_2)), $$
where $P_{H_0}$, $P_{H_1}$ are the orthogonal projections onto the corresponding
subspaces, $P_{H_{\rm res}}$ is the orthogonal projection onto $H_{\rm res}$
$$
P_{H_{\rm res}}=I_{H}-\mathcal G^*\mathcal G,
$$
where $\mathcal G$ is defined in \eqref{4.110}.
Then $\Gamma (N)$ is a cyclic subspace for the operator $U$ and the corresponding spectral function $a_\omega (\zeta):\begin{bmatrix}N_2\cr N_1\cr\end{bmatrix}
\to\begin{bmatrix}N_2\cr N_1\cr\end{bmatrix}$ is given by the formula
$$ a_\omega (\zeta)=\frac{1}{2}\begin{bmatrix}0 && -\omega (0)\cr\cr \omega (0)^* && 0 \cr
\end{bmatrix}+ {\buildrel\circ \over a}_\omega (\zeta)\qquad\qquad\qquad\qquad\qquad\qquad\qquad  \eqno(5.3) $$
$$\qquad\qquad - \frac{1}{2}\int\limits_{\mathbb T }{t+\zeta\over t-\zeta}\begin{bmatrix}\psi^* && \omega\varphi
\cr\cr \omega^*\psi^* && \varphi\cr\end{bmatrix}\ \begin{bmatrix}{\bf 1}_{E_2} && w
\cr\cr w^* && {\bf 1}_{E_1}\cr\end{bmatrix}^{[-1]}\begin{bmatrix}\psi && \psi
\omega \cr\cr \varphi^*\omega^* && \varphi^* \cr\end{bmatrix} m(dt), $$
where
$$ {\buildrel\circ\over a}_\omega=\frac{1}{2}\ {{\bf 1}_{N_2\oplus N_1}+\begin{bmatrix}
0 && \omega \cr\cr s && 0 \cr\end{bmatrix}\over {\bf 1}_{N_2\oplus N_1}-
\begin{bmatrix}0 && \omega \cr\cr s && 0 \cr\end{bmatrix}}\ = \
\begin{bmatrix}{\buildrel\circ\over \psi}-\frac{1}{2} &&
\omega{\buildrel\circ\over\varphi} \cr\cr s{\buildrel\circ
\over \psi} && {\buildrel\circ\over\varphi}- \frac{1}{2} \cr\end{bmatrix}\ , $$
$${\buildrel\circ\over\varphi}=({\bf 1}_{N_1}-s\omega)^{-1},\quad
{\buildrel\circ\over \psi} = ({\bf 1}_{N_2}-\omega s)^{-1},\eqno(5.4)$$
$$\varphi={\buildrel\circ\over\varphi}s_1=({\bf 1}_{N_1}-s\omega)^{-1}s_1,\quad
\psi=s_2{\buildrel\circ\over\psi} = s_2({\bf 1}_{N_2}-\omega s)^{-1},$$
$w$ is the characteristic function of the colligation $\mathcal A$
$$ w=s_0+s_2\omega({\bf 1}-s\omega)^{-1}s_1, \eqno(5.5) $$
$m(dt)$ is the normalized Lebesgue measure on the unit circle $\mathbb T $.
\end{theorem}

\medskip\noindent
Note that the spectral function $a_\omega$ of the
feedback coupling depends on the characteristic functions of
the coupled colligations only.

\medskip\noindent
The real part of the spectral function $a_\omega(\zeta)$ can be expressed as
$$ a_\omega(\zeta)+a_\omega(\zeta)^* = {\buildrel\circ\over a}_\omega(\zeta)+{\buildrel
\circ\over a}_\omega(\zeta)^* 
\qquad\qquad\qquad\qquad\qquad\qquad\qquad\qquad\qquad \eqno(5.6) $$

\medskip\noindent
$$ - \int\limits_{\mathbb T }{1-|\zeta|^2\over | t-\zeta|^2}\ \begin{bmatrix}\psi^*
&& \omega\varphi \cr\cr \omega^*\psi^* && \varphi \cr\end{bmatrix}\ \begin{bmatrix}{\bf 1}_{
E_2} && w \cr\cr w^* && {\bf 1}_{E_1} \cr\end{bmatrix}^{[-1]}\begin{bmatrix}
\psi && \psi\omega \cr\cr \varphi^*\omega^* && \varphi^* \cr\end{bmatrix} m(dt). $$
We will also need the following re-expression of ${\buildrel\circ\over
a}_\omega$
$$ {\buildrel\circ\over a_\omega}=\begin{bmatrix}\frac{1}{2} & {{\bf 1}_{N_2}+\omega s
\over {\bf 1}_{N_1}-\omega s} &&& \ &  \omega{\buildrel\circ\over\varphi} \cr\cr  &
s{\buildrel \circ\over\psi} &&& \frac{1}{2} & {{\bf 1}_{N_1}+s\omega\over {\bf 1}_{N_1}-
s\omega}\cr\end{bmatrix}  \ \eqno(5.7) $$
with the same notations as in (5.4).
Since $${\buildrel\circ\over a}_w + {\buildrel\circ\over a}_w^*\geq 0,$$
there exists an operator measure ${\buildrel\circ\over\sigma }(dt):N\to N$
such that
$$ {\buildrel\circ\over a}_\omega(\zeta)=\frac{1}{2}\begin{bmatrix}0 && \omega(0) \cr\cr
-\omega(0)^* && 0 \cr\end{bmatrix} + \frac{1}{2}\int\limits_{\mathbb T }{t+\zeta\over t-\zeta}\ {
\buildrel\circ\over\sigma }(dt). \ \eqno(5.8) $$

\medskip\noindent
\subsection{Property $\Vert Fx\Vert^2=D(x,x)$ yields Boundary Properties
of the Coefficient  Matrix $S$ and the Parameter $\omega$.}

\medskip\noindent
In this section we show how formula (5.6) for the spectral
function of the residual part yields boundary
properties of the coefficient matrix-function
$$S=\begin{bmatrix}s_0 && s_2 \cr\cr s_1 && s \cr\end{bmatrix}.$$
Let $w$ be a solution of the AIP and let $F^\omega:X\to H^w$ (we label it
with parameter $\omega$ in view of formula \eqref{4.110})
be the corresponding mapping. Assume that the equality
$$ \Vert F^\omega x\Vert^2_{H^w} = D(x,x) \eqno(5.9) $$
holds for all $x\in X$. By Proposition \ref{5.300},
this means that the corresponding extension $\mathcal A$ is simple, i.e., the residual
part is trivial. But then the spectral function of the residual part must
be equal to zero. Therefore, formula (5.6) along with assumption (5.9) yield
$$ {\buildrel\circ\over a}_\omega(\zeta)+{\buildrel\circ\over a}_\omega(\zeta)^*
\qquad\qquad\qquad\qquad\qquad\qquad\qquad\qquad\qquad\qquad\qquad \eqno(5.10)$$

$$ =
\int\limits_{\mathbb T }{1-|\zeta|^2\over | t-\zeta|^2}\ \begin{bmatrix}\psi^* &&
\omega\varphi \cr\cr \omega^*\psi^* && \varphi\cr\end{bmatrix}\ \begin{bmatrix}{\bf 1}_{E_2}
&& w \cr\cr \omega^*\psi^* && \varphi \cr\end{bmatrix}^{[-1]} \begin{bmatrix}\psi
&& \psi\omega \cr\cr \varphi^*\omega^* && \varphi^*\cr\end{bmatrix}m(dt). $$
The latter is equivalent to the following two properties:
\begin{enumerate}
\item\ ${\buildrel\circ\over\sigma }_\omega(dt)$ is absolutely continuous with
respect to the Lebesgue measure;
\item\ ${\buildrel\circ\over a}_w(t)+{\buildrel\circ\over a}_\omega(t)^* =
\begin{bmatrix}\psi^* && \omega\varphi \cr\cr \omega^*\psi^* && \varphi\cr\end{bmatrix}\
\begin{bmatrix}{\bf 1}_{E_2} && w\cr\cr w^* && {\bf 1}_{E_1}\cr\end{bmatrix}^{[
-1]}\begin{bmatrix}\psi && \psi\omega \cr\cr \varphi^*\omega^* && \varphi^* \cr
\end{bmatrix}$,
\end{enumerate}
almost everywhere on $\mathbb T $. Thus, these two properties are equivalent to
$$ \Vert F^\omega x\Vert^2_{H^w}=D(x,x), \ \forall\ x\in X. $$
%
Problem \ref{2.10} and Problem \ref{2.350} possess the property
$$ \Vert F^\omega x\Vert^2_{H^w}=D(x,x), \ \forall\ x\in X $$
for every solution $w$ and every parameter $\omega$ that produces this $w$ via formula (5.5)
(actually, for these examples formula (5.5) gives a one to one correspondence between $\omega$
and $w$). Therefore, properties 1 and 2 above hold for every parameter $\omega$ in those examples.

\medskip\noindent
Property 2 itself is equivalent to vanishing of the absolutely continuous
part of the measure $\sigma _\omega(dt)$ that corresponds to $a_\omega(\zeta)$.
It is the case in
Problem \ref{2.1000} that ${\buildrel\circ\over\sigma }_\omega(dt)$
can be supported
by a single point for every parameter $\omega$. Hence, $\sigma _\omega(dt)$ has trivial
absolutely continuous part for every parameter $\omega$ in this example.
Thus, we have Property 2 for every parameter $\omega$ in 
Problem \ref{2.1000}. Property~1
holds for some parameters and does not hold for the others in this
problem. Analysis shows that in this problem the residual part
is nontrivial if and only if
$$
\underset{\zeta\to t_0}{\lim} \omega(\zeta)=\overline{s(t_0)},
$$
where $s$ is the right-bottom entry of the coefficient matrix $S$,
and
$$
\underset{\zeta\to t_0}{\lim} \frac {1-|\omega(\zeta)|^2}{1-|\zeta|^2} <\infty .
$$
Note that $|s(t_0)|=1$ always in this problem.

\medskip\noindent
We reformulate the above Properties 1 and 2 using
formulas (5.6) and (5.7) of this Lecture. Property 1 is equivalent
to the following property

\medskip\noindent
{\bf 1'.} Measures ${\buildrel\circ\over\sigma }_w^1(dt)$ and ${\buildrel\circ\over\sigma }
_\omega^2(dt)$ are absolutely continuous with respect to the Lebesgue measures,
where
$$ {{\bf 1}_{N_2}+\omega(\zeta)s(\zeta)\over {\bf 1}_{N_2}-\omega(\zeta)s(\zeta)}\ =\ \int_
{\mathbb T }{t+\zeta\over t-\zeta}\ {\buildrel\circ\over\sigma }^2_\omega(dt) $$
and
$$ {{\bf 1}_{N_1}+s(\zeta)\omega(\zeta)\over {\bf 1}_{N_1}-s(\zeta)\omega(\zeta)}\ =\ \int_
{\mathbb T }{t+\zeta\over t-\zeta}\ {\buildrel\circ\over\sigma }^2_\omega(dt).$$

\medskip\noindent
A special case of Property 2, when
parameter $\omega=0$, reads as follows

\medskip\noindent
{\bf 2'.} $\begin{bmatrix}{\bf 1}_{N_2} && s^* \cr\cr s && {\bf 1}_{N_1} \cr
\end{bmatrix}\ = \ \begin{bmatrix}s^*_2 && 0 \cr\cr  0 && s_1 \cr\end{bmatrix} \
\begin{bmatrix}{\bf 1}_{E_2} && s_0 \cr\cr s_0^* && {\bf 1}_{E_1} \cr\end{bmatrix}
^{[-1]}\ \begin{bmatrix}s_2 && 0 \cr\cr  0 && s_1^* \cr\end{bmatrix}$

\medskip\noindent
almost everywhere on $\mathbb T $.

\medskip\noindent
It was shown in \cite{9,10} (under the assumptions that $\dim E_1<\infty$ and
$\dim E_2 < \infty$) that converse is also true, i.e. Property {\bf 2'} implies {\bf 2} for every
parameter $\omega$. It was also shown in \cite{9,10} under the same assumptions ($
\dim E_1<\infty$ and $\dim E_2<\infty$) that Property {\bf 2'} is in turn equivalent
to this property

\medskip\noindent
{\bf 2".}\ \ ${\rm rank}\ ({\bf 1}_{E_1,\oplus N_2}-S^*S)={\rm rank}\ ({\bf 1}_{E_1}- s^*_0s_0)
-\dim N_1$,\quad a.e. on $\mathbb T $.

\medskip\noindent
This result contains, in particular , the classical Nevanlinna-Adamjan-Arov-Krein
theorem:\ \ if $\dim N_1=\dim E_1$ (i.e., the problem is completely
indeterminate), then
$$ {\bf 1} - S^* S = 0 $$
a.e. on $\mathbb T $. That is, the coefficient matrix $S$ in inner.

\section{Lecture 6: Properties of the Coefficient Matrices via Dense Sets in the Function Model Spaces.}

\medskip\noindent
Let $S$ be the characteristic function \eqref{4.75} of the universal colligation $\mathcal A_0$ \eqref{4.60}
$$ S(\zeta)=\begin{bmatrix}s_0(\zeta) && s_2(\zeta)\cr\cr s_1(\zeta) && s(\zeta)
\cr\end{bmatrix}:\begin{bmatrix}E_1 \cr\cr N_2 \cr\end{bmatrix}\to\begin{bmatrix}
E_2 \cr\cr N_1 \cr\end{bmatrix}. $$
It serves as the coefficient matrix in the parametrization formula \eqref{4.140}.
Since $S$ is a Schur class operator-function, there exists de Branges-Rovnyak function space $H^S$
associated to it (similar to definition in Section 1.2).
We discuss here an approach that employs some dense sets in the space $H^S$.
This approach was applied to Sarason problem in \cite{KhReg}, to Nehari problem in \cite{KhExposed, KhNeh} and to the general Commutant Lifting problem in \cite{BallK1, BallK2}.

\subsection{The data of AIP suggest a dense set in $H^S$.}

\medskip\noindent
Let ${\mathcal G}^0$ be the Fourier representation of the colligation
$\mathcal A_0$ defined in \eqref{4.80}.
${\mathcal G}^0$ maps the
space $H^0$ onto the de Branges-Rovnyak function space $H^S$ contractively.
We define mapping $F^S:X\to H^S$ as
$$ F^S x \ {\buildrel def\over =}\ {\mathcal G}^0[x], \ \eqno(6.1) $$
where $[x]$ is defined in \eqref{3.5}, \eqref{3.7}.
Observe that the lineal $\{F^Sx, \ x\in X\}$ is
dense in $H^S$, since the lineal $\{[x], \
x\in X\}$ is dense in $H^0$.

\medskip\noindent
It follows from definition \eqref{4.10}-\eqref{4.30} of the colligation $\mathcal A_0$
and definition \eqref{4.80} of the Fourier representation ${\mathcal G}_0$
that $F^S$ possess the properties similar to \eqref{1.30} and \eqref{1.40}:

\medskip
i) \hskip 2cm $\Vert F^Sx\Vert^2_{H^S}\leq D(x,x), \ \ x\in X$;

ii)\hskip 2cm $tF^ST_2x-F^ST_1x= \ \begin{bmatrix}{\bf 1}_{E_2\oplus N_1}
& & S \cr\cr S^* && {\bf 1}_{E_1\oplus N_2}\cr\end{bmatrix}\ \begin{bmatrix}
-M_2x \cr\cr 0 \cr\cr M_1x \cr\cr 0 \cr\end{bmatrix}$ ,

\medskip\noindent
a.e. on $\mathbb T $, where the first zero stands for the zero vector of the
space $N_1$ and the second zero stands for the zero vector of the space
$N_2$.

\medskip\noindent
Under the Special Case assumptions \eqref{1.50},
property ii) can be re-expressed as follows:
$$ (F^S_+x)(\zeta)=\left(S(\zeta)\begin{bmatrix}M_1\cr\cr 0 \cr
\end{bmatrix} - \begin{bmatrix}M_2 \cr\cr 0 \cr\end{bmatrix}\right)\
(\zeta T_2-T_1)^{-1}x $$
$$ (F^S_-x)(\zeta)=\bar\zeta\left(\begin{bmatrix}M_1 \cr\cr 0 \cr\end{bmatrix}-S
(\zeta)^*\begin{bmatrix}M_2\cr\cr 0 \cr\end{bmatrix}\right)(T_2-\bar\zeta T_1
)^{-1}x,\quad |\zeta|<1. \ \eqno(6.2) $$
Here
$$\begin{bmatrix}M_1 \cr\cr 0 \cr\end{bmatrix}:X\to
\begin{bmatrix}E_1 \cr\cr N_2 \cr\end{bmatrix},\qquad \begin{bmatrix}M_2 \cr\cr
0 \cr\end{bmatrix} :X\to \begin{bmatrix}E_2 \cr\cr N_1 \cr\end{bmatrix}.$$

\medskip\noindent
The structure of the dense set $\{F^Sx,x\in X\}$ in $H^S$ carries all the
information about properties of the matrix-function $S$.
It will be demonstrated in the next sections how this approach works for Sarason problem.

\subsection{Coefficient Matrix of the Sarason Problem and Associated Function Model Space.}

\medskip\noindent
In Section \ref{Example 2} we considered the scalar-valued case of the Sarason problem.
We consider more general case now.
Let $\theta$ be an inner operator-function, $\theta(\zeta):E^\prime_2\to E_2, \ \
\theta^*\theta={\bf 1}_{E_2^\prime}$, a.e. on $\mathbb T ,\ E^\prime_2$ and $E_2$ be
separable Hilbert spaces. Let $K_\theta=H^2_+(E_2)\ominus\theta H^2_+(E^\prime
_2)$, where $H^2_+(E)$ stands for the vector Hardy space with coefficients
in the space $E$. Let $T^*_\theta x=P_+\bar t x,\ x\in K_\theta$. Let $W^*$ be
a contractive operator, $W^*:K_\theta\to H^2_+(E_1)$, where $E_1$ is a
separable Hilbert space, such that $W^*T_\theta^*=P_+\bar t W^*$. Consider
the following interpolation problem: find all
the Schur class functions $w(\zeta):E_1\to E_2$ that satisfy
$$ W^*x=P_+w^*x, \ x\in K_\theta\ . $$
The associated AIP data are:
$$ X=K_\theta,\ T_1=I_X,\ T_2=T_\theta^*, $$
$E_1$ and $E_2$ are the spaces introduced above, $M_1x=(W^*x)(0),\ M_2x
=x(0)$, where the latter notations stands for the value at $0$ of an
$H^2_+(E_1)$ and an $H^2_+(E_2)$ function, respectively,
$$
D(x,x)=\langle (I-W W^*)x, x \rangle_{_{K_\theta}}.
$$
The Fourier representation $F^w$ is unique for every solution $w$ and can be expressed as
$$ F^wx = \begin{bmatrix}{\bf 1}_{E_2} & & w \cr\cr w^* && {\bf 1}_{E_1}\cr
\end{bmatrix}\ \begin{bmatrix}x \cr\cr -W^*x \cr\end{bmatrix}. $$
The Fourier representation $F^S$ of the colligation $\mathcal A_0$ (see (6.2))
can also be expressed as
$$ F^S x=\begin{bmatrix}{\bf 1}_{E_2\oplus N_1} && S \cr\cr S^* && {\bf 1}
_{E_1\oplus N_2}\cr\end{bmatrix}\ \begin{bmatrix}x \cr\cr 0 \cr\cr -W^*x
\cr\cr  0 \cr\end{bmatrix},\ \ x\in K_\theta\ . $$
Since $s_0$ is also a solution (corresponding to the parameter $\omega=0$), we can
write the latter expression as
$$ F^S x=\begin{bmatrix}{\bf 1}_{E_2\oplus N_1} && S \cr\cr S^* && {\bf 1}_
{E_1\oplus N_2}\cr\end{bmatrix}\ \begin{bmatrix}x \cr\cr 0 \cr\cr -P_+s^*_0x
\cr\cr 0 \cr\end{bmatrix},\ x\in K_\theta\ . \ \eqno(6.3) $$
As it was observed in Section 6.1, the set
$$ \{F^Sx,\ x\in K_\theta\}\subseteq H^S $$
is dense in $H^S$.

\subsection{Properties of the Coefficient Matrices of the Sarason
Problem.}

\medskip\noindent
Observe first, that the bottom entry of the vector $F^Sx$ (i.e. the
second component of the vector $F^S_-x$ ) in (6.3) is $s^*_2x$, and it is an
$H^2_-(N_2)$ function, since $F^Sx\in H^S$.
Thus,
$$ s^*_2x\in H^2_-(N_2),\quad \forall  x\in K_\theta . $$ 
In other words,
$$ \langle s_2^*x,\ h_+\rangle_{L^2(N_2)}=0 $$
for all $h_+\in H^2_+(N_2)$ and $x\in K_\theta$. Equivalently,
$$ \langle x,\ s_2h_+\rangle_{L^2(E_2)}=0 $$
for all $h_+\in H^2_+(N_2)$ and $x\in K_\theta$. Hence, $s_2h_+\in\theta H^2_+
(E_2^\prime)$ for all $h_+\in H^2_+(N_2)$. This means that
$$ s_2=\theta\tilde s_2 $$
for a Schur class function $\tilde s_2:N_2\to E_2^\prime$.
The next theorem shows how the denseness property works.

\begin{theorem}[\cite{1}]
Let the $S=\begin{bmatrix}s_0 && \theta
\tilde s_2 \cr\cr s_1 && s \cr\end{bmatrix}$ be the coefficient matrix
of the Sarason problem. Then
$s_1$ is an  outer function $($i.e. the lineal $s_1H^2_+(E_1)$ is dense
in $H^2_+(N_1))$ and $\tilde s_2$ is a $*$-outer function $($i.e. the
lineal $\tilde s_2^*H^2_-(E^\prime_2)$ is dense in $H^2_-(N_2))$.
\end{theorem}

\medskip\noindent
{\bf Sketch of the proof.} To prove that $s_1$ is an outer
function  we need to check that the assumption
$$ P_+s_1^*h_+=0,\ \ h_+\in H^2_+(N_1) \ \eqno(6.4) $$
implies $h_+=0$. Consider the vector
$$ \begin{bmatrix}
{\bf 1} && S \cr\cr S^* && {\bf 1}
\end{bmatrix}\
\begin{bmatrix}
\begin{bmatrix}0 \\ h_+\end{bmatrix}\\ \\
-P_+S^*\begin{bmatrix}0 \\ h_+ \end{bmatrix}
\end{bmatrix}
=\begin{bmatrix}{\bf 1} && S \cr\cr S^* && {\bf 1}\cr
\end{bmatrix}\ \begin{bmatrix}0 \cr\cr h_+ \cr\cr -P_+s^*_1h_+ \cr\cr -P_+
s^*h_+ \cr
\end{bmatrix}. \
\eqno(6.5) $$
It belongs to $H^S$ for every $h_+\in H^2_+(N_1)$. Assumption
(6.4) makes it orthogonal to the lineal (6.3). Since (6.3) is dense in $H^S$,
the orthogonality forces (6.5) to be zero. The latter in turn yields
$h_+=0$.

\medskip\noindent
We return now to the Problem \ref{2.350} 
i.e., to the scalar Sarason problem $(E_1 =E_2=E_2^\prime=\mathbb C^1)$. Assume also
that the problem is indeterminate (i.e. permits a non-unique solution,
i.e. $N_1=N_2=\mathbb C^1$).  In this case $S$ is an inner matrix-function (see Lecture 5) and it can be normalized so that
$s_1=\tilde s_2=a$ is an outer function.
The following criterion is a consequence of the denseness of the lineal (6.3) in $H^S$.

\begin{theorem}[\cite{KhReg, KhExposed}, also \cite{BallK2}] A $2\times 2$ inner matrix function
$$ S(\zeta)=\begin{bmatrix}s_0 && s_2 \cr\cr s_1 && s \cr\end{bmatrix}\ =\
\begin{bmatrix}s_0 && \theta a \cr\cr a && s \cr\end{bmatrix}, $$
$($where $a$ is an outer function, $s(0)=0)$ is the coefficient matrix of an indeterminate Sarason
Problem \ref{2.350}
if and only if
$$ \begin{bmatrix}P_-\bar s_1 \cr\cr \bar s\cr\end{bmatrix}\in {\rm clos}\left\{
\begin{bmatrix}P_-\bar s_0x \cr\cr \bar s_2 x \cr\end{bmatrix},\ x\in K_\theta
\right\}, \eqno (6.6)$$
where the closure is understood in the $L^2$ sense.
\end{theorem}
\begin{remark}{\rm
Property (6.6) is equivalent to the following one
$$ \inf_{x\in K_\theta}\Vert P_- S^*\left(\begin{bmatrix}0 \cr\cr 1 \cr\end{bmatrix}
\ - \begin{bmatrix}x \cr\cr 0 \cr\end{bmatrix}\right) \Vert^2_{H^2_-}=0 , $$
equivalently,
$$ \inf_{x\in K_\theta}\left\langle  SP_- S^*\left(\begin{bmatrix}0 \cr\cr 1 \cr\end{bmatrix}
\ - \begin{bmatrix}x \cr\cr 0 \cr\end{bmatrix}\right), \ \left( \begin{bmatrix}0 \cr\cr 1 \cr\end{bmatrix} -
\begin{bmatrix}x \cr\cr 0 \cr\end{bmatrix}\right) \
\right\rangle _{L^2}=0. \ \eqno(6.7) $$
Further simplification of property (6.7) is unknown. A discussion on this matter was given in \cite{NATO}.
An application to uniqueness of the inverse scattering for CMV matrices was given in \cite{VoYu, GKY}.}
\end{remark}


\begin{thebibliography}{99}

\bibitem{23}
V.M. Adamjan, D.Z. Arov, M.G. Krein,
\textit{Infinite Hankel matrices and generalized problem of
Carath\'{e}odory-Fej\'{e}r and I. Schur.},
Funkt. Anal. i Prilozhen.   {\bf 2:4} (1968) 1-17, Russian.
English Transl.,   Funct. Analysis and Appl.  {\bf 2} (1968) 269-281.

\bibitem{24}
V.M. Adamjan, D.Z. Arov, M.G. Krein,
\textit{Some function theoretic problems connected with the theory of
spectral measures of isometric operators}, in: Lecture Notes in
Mathematics, Springer-Verlag, {\bf 1043} (1984) 160-163; or
{\bf 1573} (1994) 183-185.

\bibitem{ArGr}
D. Z. Arov and L. Z. Grossman, \textit{Scattering matrices in the
theory of extensions of isometric operators}, Dokl. Akad. Nauk SSSR 270:1 (1983),
17-20. Russian; English translation in Soviet Math. Dokl. 27:3 (1983), 518-521.

\bibitem{Ball1978}
 J. A. Ball, \textit{Factorization and model theory for contraction operators with
unitary part}, Mem. Amer. Math. Soc. {\bf 13} (1978), no. 198

\bibitem{BallK1}
J. Ball, A. Kheifets,
\textit{The Inverse Commutant Lifting Problem.
I: Coordinate-Free Formalism}, Integral Equations and Operator Theory {\bf 70} (2011) 17-62.

\bibitem{BallK2}
J. Ball, A. Kheifets,
\textit{The Inverse Commutant Lifting Problem. II:
Hellinger Functional-Model Spaces}, Complex Analysis and Operator Theory {\bf 7} (2013) 873-907.

\bibitem{BoKh1}
V. Bolotnikov, A. Kheifets,
\textit{A higher order analogue of the Caratheodory-Julia theorem},
Journal of Functional Analysis {\bf 237} 1 (2006) 350--371

\bibitem{BoKh2}
V. Bolotnikov, A. Kheifets,
\textit{A higher order Caratheodory-Julia theorem and related boundary interpolation problem},
Operator Theory: Advances and Applications {\bf 179} (2008) 63--102

\bibitem{BrLiv}
M. S. Brodskii, M. S. Liv\v sic, \textit{Spectral analysis of non-selfadjoint operators and intermediate systems},
(Russian) Uspehi Mat. Nauk (N.S.) {\bf 13} (1958) no. 1(79), 3-85.
Engl. Transl.: Amer. Math. Soc. Transl. (2) {\bf 13} (1960) 265-346

\bibitem{deBr}
L. de Branges, \textit{Factorization and invariant subspaces}, J. Math. Anal. Appl. {\bf 29} (1970) 163-200

\bibitem{GKY}
L. Golinskii, A. Kheifets, P. Yuditskii,
\textit{Scattering theory for CMV matrices: uniqueness, Helson-Szeg\H o and strong Szeg\H o theorems}, Integral Equations Operator Theory
{\bf 69} (2011) no. 4 479-508.

\bibitem{KKY}
V. E. Katsnel'son, A. Y. Kheifets, and P. M. Yuditskii,
\textit{Abstract Interpolation Problem and Theory of Extensions of Isometric Operators},
in Operators in function spaces and problems in function theory
{\bf 146} (1987) 83-96 Naukova Dumka, Kiev (Russian); English translation
in Oper. Theory Adv. Appl. {\bf 95} (1997) 283-298 Birkh\" auser, Basel.

\bibitem{6}
A.Ya. Kheifets, \textit{Parseval equality in abstract
interpolation problem and coupling of open systems}, Teor. Funk.,
Funk. Anal. i ikh Prilozhen {\bf 49} (1988) 112--120, {\bf 50}
(1988) 98--103, Russian. English transl.:   J. Sov. Math.  {\bf
49}, 4 (1990) 1114--1120, {\bf 49}, 6 (1990) 1307--1310.

\bibitem{1}
A.Ya. Kheifets, \textit{Generalized bitangential
Schur - Nevanlinna - Pick problem, related Parseval equality and scattering
operator}, deposited in VINITI, 11.05.1989, No. 3108--B89 Dep.,
1--60, 1989, Russian.

\bibitem{9}
A. Ya. Kheifets, \textit{Scattering matrices and Parseval
equality in Abstract Interpolation Problem}, PhD thesis, 1990, Russian.

\bibitem{10}
A.Ya. Kheifets, \textit{Nevanlinna - Adamjan - Arov - Krein
theorem in semi - determinate case},   Teor. Funkt., Funkt. Anal. i
ikh Prilozhen.  {\bf 56} (1991) 128--137, Russian. Engl. transl.:
Journal of Mathematical Sciences {\bf 76}, 4 (1995) 2542--2549.


\bibitem{KhReg}
A.Ya. Kheifets, \textit{Regularization of gamma-generating pairs and exposed points in $H^1$}, Journal of Functional Analysis {\bf 130} 2 (1995) 310--333.

\bibitem{KhHamb}
A.Ya. Kheifets, \textit{Hamburger moment problem: Parseval equality and Arov-singularity}, Journal of Functional Analysis, {\bf 141} 2 (1996) 374--420

\bibitem{KhExposed}
A. Kheifets, \textit{Nehari problem and exposed points of the unit ball in the Hardy space} , Proceedings of the Ashquelon conference, Israel, 1996, in Israel Mathematical Conference Proceedings, {\bf 11} (1997) 145-151


\bibitem{KhBer}
A. Kheifets, \textit{Abstract interpolation problem and some applications}, Lecture notes given in framework of Holomorphic Spaces semester, fall 1995, in: Holomorphic Spaces (S. Axler, J. McCarthy, D. Sarason editors), MSRI Publications, {\bf 33} (1998) 351-381, Berkeley, California.


\bibitem{KhNeh}
A. Kheifets, \textit{Parameterization of solutions of the Nehari Problem and nonorthogonal dynamics}, Operator Theory: Advances and Applications, {\bf 115}  (2000) 213-233, Birkhauser Verlag, Basel.

\bibitem{NATO}
A. Kheifets,
\textit{Density of domain of the weighted Hilbert transform},
Byrnes, James S. (ed.), Twentieth century harmonic analysis - a celebration. Proceedings of the NATO Advanced Study Institute, Il Ciocco, Italy, July 2-15, 2000. Dordrecht: Kluwer Academic Publishers. NATO Sci. Ser. II, Math. Phys. Chem. {\bf 33}, 374-375 (2001). Zentrallblatt Review Zbl 0991.44501, https://zbmath.org/?q=an:0991.44501

\bibitem{KhY-0}
A. Ya. Kheifets, P.M. Yuditskii,
\textit{Interpolation of operator that commutes with compressed shift by Schur class functions},
Teor. Funk., Funk. Anal. i ikh Prilozhen, {\bf 40} (1983) 129--136, Russian

\bibitem{KhY}
A. Y. Kheifets and P. M. Yuditskii, \textit{ Analysis and
extension of V. P. Potapov's approach to interpolation problems with applications
to the generalized bi-tangential Schur-Nevanlinna-Pick problem and J-inner-outer
factorization}, Oper. Theory Adv. Appl. {\bf 72} (1994) 133-161.

\bibitem{VoYu}
 A. Volberg, P. Yuditskii, \textit{On the inverse scattering problem for Jacobi matrices with the spectrum on an interval, a finite system of intervals or a Cantor set of positive length}, Comm. Math. Phys. {\bf 226} (2002) no. 3 567-605

\bibitem{Sarason}
D. Sarason, Sub-Hardy Hilbert spaces in the unit disk, University of
Arkansas Lecture Notes in the Mathematical Sciences, Wiley, New York, 1994.


\bibitem{NF}
 B. Sz.-Nagy, C. Foia\c s, \textit{Harmonic analysis of operators on Hilbert space}, Translated from the French and revised North-Holland Publishing Co., Amsterdam-London; American Elsevier Publishing Co., Inc., New York; Akademiai Kiado, Budapest 1970 xiii+389 pp.

\bibitem{22}
Yu. L. \v Smulian, \textit{Some points of the theory
of operators with finite nonhermitian rank},
  Matem. Sborn.  {\bf 57} (99) (1962) 105--136, Russian.

\bibitem{Schv}
Ja. S. \v Svarcman, On invariant subspaces of a dissipative operator and the divisors of its characteristic function, Funkcional. Anal, i Priložen. {\bf 4} (1970), no. 4, 85-86. Engl. transl.: Functional Anal. Appl. {\bf 4} (1970), 342-343.

\end{thebibliography}
\end{document}